\documentclass{pomona2004}
\def\ppages{1--13}

\usepackage{graphicx}
\usepackage{amscd}

\theoremstyle{plain}
\newtheorem{theorem}{Theorem}

\newtheorem{lemma}{Lemma}

\theoremstyle{definition}

\newtheorem{remark}{Remark}

\input{epsf}

\title[Colored Coalescent Theory]{Colored Coalescent Theory}
\email{tianjj@mbi.ohio-state.edu}
\email{xl@math.ucr.edu}
\author[Jianjun Tian and  Xiao-Song Lin]{}

\begin{document}
\maketitle

\centerline{\scshape Jianjun Tian}
{\footnotesize \centerline{}
  \centerline{Mathematical Biociences Institute}
   \centerline{The Ohio State University}
   \centerline{Columbus, OH 43210, USA}
} 

\bigskip
\centerline{\scshape Xiao-Song Lin}
{\footnotesize \centerline{}
  \centerline{Department of Mathematics}
   \centerline{University of California, Riverside}
   \centerline{Riverside, CA 92521, USA}
} 



 


\begin{quote}{\normalfont\fontsize{8}{10}\selectfont
{Abstract.} 
We introduce a colored coalescent process which recovers random colored genealogical trees. Here a colored genealogical tree has its vertices colored black or white. Moving backward along the colored genealogical tree, the color of vertices may change only when two vertices coalesce. Explicit computations of the expectation and the cumulative distribution function of the coalescent time are carried out. For example, when $x=1/2$, for a sample of $n$ colored individuals, the expected time for the colored coalescent process to reach a black MRAC or a white MRAC, respectively, is $3-2/n$. On the other hand, the expected time for the colored coalescent process to reach a MRAC, either black or white, is $2-2/n$, which is the same as that for the standard Kingman coalescent process.
\par}
\end{quote}

\section{Introduction}

In the last twenty years, coalescent theory has been developed into a powerful analytical tool for population genetics. This theory is especially significant with the rapid accumulation of DNA sequence data. First formulated in the seminal work of Kingman in 1982 \cite{King1,King2}, coalescent theory offers various sample-based and highly efficient statistical methods for analyzing molecular data such as DNA sequence samples. For recent reviews as well as extensive references of coalescent theory, see \cite{FL} and \cite{RN}. A nice introduction to coalescent theory can be found in \cite{no}.
\smallskip

Mathematically, coalescent theory studies stochastic processes leading to the most recent common ancestor (MRCA) from a sample under various coalescent models. If one thinks of the more commonly studied branching processes as 
stochastic models of generating random trees from their roots, coalescent processes can be thought of as the inverse processes which recover random trees from their leaves. In a more elaborated version crucial for population genetics, a coalescent process is usually superimposed with a mutation process.
This mutation process can be thought of as an independent Poisson process running on the random tree generated by the coalescent process, with the edge lengths of the random tree serving as the time scale for the mutation process.
\smallskip

In this article, we introduce a coalescent process which generated random colored trees. Here a {\it coloring} of a tree is to color the vertices of the tree by two colors, black ($B$) and white ($W$), such that if two vertices are joint by an edge, they may have different colors only when the vertex closer to the root is a branching point of the tree. Thus, in recovering a random colored tree using this coalescent process, we may end up at a colored tree with the root colored black or colored white. The quantities which we are interested in include the following: the probabilities for the coalescent process to reach a black or white root, respectively; the mean and the cumulative distribution function of the coalescent time, which is the time elapsed before the coalescent process reaches a black root or a white root, respectively; etc. 

\smallskip
Our goal in the mathematical development of colored coalescent theory is to understand the geographical models of the origin of the human being.  Let us take two sets of gene sequences from the different geographical locations.  Designate each by different color.  We then mix them as a population sample.  Based upon the sample, we look at which location the most recent common ancestor was arisen from.  Our theory can give a answer theoretically.  We hope to incorporate genetic data analysis of gene genealogies into the colored coalescent framework in our future study.
\smallskip

\smallskip

The Wright-Fisher model in population genetics assumes discrete, non-overlapping generations $G_0,G_1,G_2,\dots$ in which each generation contains a fixed number $N$ of individuals. In a so-called haploid population, each member in $G_{i+1}$ is the child of exactly one member in $G_i$, but the number of children born to the $j$'s member of $G_i$ is a random variable $\nu_j$ satisfying the symmetric multinomial distribution
$$\Pr\{\nu_j=n_j\,;\,j=1,2,\dots,N\}=\frac{N!}{n_1!n_2!\cdot n_N!N^N}.$$
\smallskip

We additionally assume that each individual in a generation has two possible colors $B$ and $W$. In the next generation, if a member is the only child of its parent, then this child will inherit the color of its parent. But when a parent has more than one child in the next generation, the color of children of that common parent satisfies a binomial distribution. More specifically, for a parent with $k$ children in the next generation, $k>1$, let $b$ be the number of children with $B$ color and $w$ be the number of children with $W$ color (so that $b+w=k$), we have
\begin{equation}\label{color1}
\begin{aligned}
&\Pr\{b=k_1,w=k_2\,;\,\text{the parent has color $B$}\}
=\binom{k}{k_1}p^{k_1}(1-p)^{k_2},\\
&\Pr\{b=k_1,w=k_2\,;\,\text{the parent has color $W$}\}     
=\binom{k}{k_1}(1-q)^{k_1}q^{k_2}  
\end{aligned}
\end{equation}
where $0\leq p,q\leq1$.

Following the same argument as in \cite{King1,King2}, we have a limiting coalescent process for a sample of $n$ colored individuals when $N\rightarrow\infty$. In this limiting coalescent process, one only allow to have two individuals in the sample to coalesce. When two colored individuals coalesce, the probability of the color of their common parent can be calculated according to Equation (\ref{color1}). Assuming that we may express those probabilities of various cases in the following multiplication rules: 
\begin{equation}\label{color2}
\begin{aligned}
BB&=xB+(1-x)W \\
BW&=(1-x)B+xW \\
WW&=xB+(1-x)W,
\end{aligned}
\end{equation}
then we must have $p=1-q$ and $x=1/2$. For other values of $x$, the coalescent process itself is still defined, although we no longer have it as the backward process of a branching genealogical process. Throughout this article, we will assume that  $0<x<1$. We call this limiting coalescent process {\it the colored coalescent process}. See Figure 1 for an example of a colored genealogical relation.
\smallskip

\bigskip
\centerline{\epsfxsize=1in\epsfbox{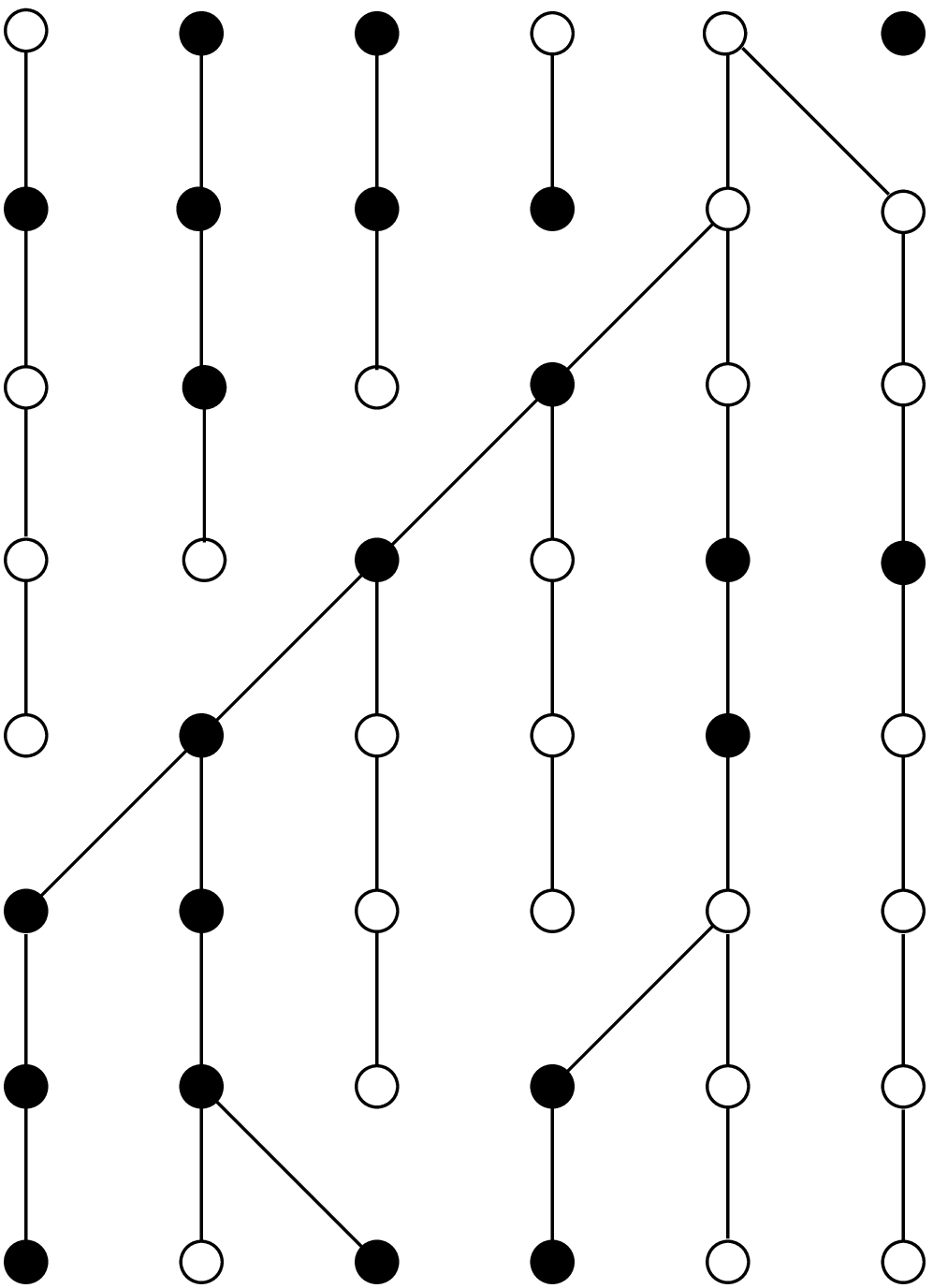}}
\medskip
\centerline{\small Figure 1. A colored genealogical relation.}
\bigskip 

The stochastic character of a state of the colored coalescent process with the parameter $x$ turns out to depend only on the number of individuals in this state colored by $B$. If we start the colored coalescent process with a sample of $n$ colored individuals, the initial state can be denoted by a pair of non-negative integers $(k,n-k)$, where $k$ is the number of individuals colored by $B$ and $n-k$ is the number of individuals colored by $W$. The absorbing states of the colored coalescent process are $(0,1)$ (a white root) and $(1,0)$ (a black root), respectively.  Our main result is Theorem 2, which states the coalescent probabilities $P_{(n_1,n_2)(0,1)}$ and $P_{(n_1,n_2)(1,0)}$ from different configurations of the sample to different functions of the colored coalescent times, $T_{(n_1,n_2)(0,1)}$ and $T_{(n_1,n_2)(1,0)}$.  For example, when $x=1/2$, the probabilities to reach $(0,1)$ and $(1,0)$, respectively, are both $1/2$.  The mean time to reach $(0,1)$ is $3-2/n$, and the mean time to reach $(1,0)$ is $3-2/n$.  Finally, according to Lemma 2, the mean time to reach either $(0,1)$ or $(1,0)$ is $2-2/n$, which is the same as in the classical coalescent theory. This shows the difference and similarity between our colored coalescent theory and the classical coalescent theory.

In Section 2, we discuss in details the colored coalescent process. The technique of lumping turns out to be very important to simplify computations involved. So we also include in Section 3 a general discussion about lumping of Markov processes. In Sections 4 and 5, we use the lumping technique to calculate various statistic parameters of the colored coalescent process. In the paper \cite{tian}, a comprehensive coalescent model which could accommodate the situation of mutations is studied.

\bigskip

\section{A colored coalescent model and some of its basic parameters}

We start with the standard Kingman coalescent process \cite{King1, King2}. So we have a sample
of $n$ individuals. In a unit time $\Delta t$, two of these $n$ individuals 
may 
coalesce with probability $\binom{n}{2}\Delta t$. Let $\Delta t\rightarrow 0$,
we arrive at a continuous stochastic process where the probability that
no coalescent event happens within the time interval $[0,t]$ is
$$e^{-\binom{n}{2} t}.$$
Here a {\it coalescent event} is the event that two of the $n$ 
individuals coalesce into one.
Thus the expectation of the time $\tau_n$ before the first coalescent 
event happens
is $E(\tau_n)={\binom{n}{2}}^{-1}.$
\smallskip

After the first coalescent event, the process will continue with a new 
generation of $n-1$ individuals. The expected time elapsed between the
first coalescent event and the second one is then 
$E(\tau_{n-1})={\binom{n-1}{2}}^{-1}$. The process continues with the number
of individuals in each generation decreases consecutively from $n$ to $1$, 
when a most recent
common ancestor (MRCA) is reached. The expected time elapsed between the 
$k$-th coalescent event and the $(k-1)$-th coalescent event is 
$E(\tau_k)=\binom{k}{2}$, for $k=n,n-1,\dots,2$. The random variable 
$$T_n=\sum_{k=2}^n\tau_k$$
is the time elapsed before the coalescent process reaches a MRCA. Its
expectation is
$$E(T_n)=\sum_{k=2}^n{\binom{k}{2}}^{-1}=2-\frac{2}{n}.$$
\smallskip

We now add an additional feature to this coalescent process by assuming that
an individual in each generation can have two colors,
black ($B$) and white ($W$).
At the initial stage of the process, each of the $n$ individuals in the 
current generation 
are given a color, $B$ or $W$. 
The process then runs as before. When two individuals coalesce, the color of
their common ancestor is determined by Equation (\ref{color2}). 
In other words, when an individual colored by $B$ and another individual
colored by $B$ coalesce, the probability that their common ancestor is
colored by $B$ is $x$ and colored by $W$ is $\overline{x}:=1-x.$ Other
situations are interpreted similarly. Furthermore, 
non-coalescent individuals will
keep their colors unchanged after an coalescent event. 
\smallskip

A typical way to study the coalescent process is to consider a 
death process on the lattice points $n,n-1,\dots,2,1$ on a line, with the 
time spending on the point $k$ exponentially distributed with the mean 
${\binom{k}{2}}^{-1}$. We will do the same for our colored coalescent 
process.
\smallskip

Fix an integer $n>0$. Consider a proper subset $\mathcal L$ of the plane 
integer lattice:
$$\mathcal L=\{(k,l)\in\mathbb Z\times\mathbb Z\,;\,k\geq 0, l\geq 0,
0<k+l\leq n\}.$$ 
A point $(k,l)\in\mathcal L$ represents a colored generation of
$k+l$ individuals, with the number of $B$-colored individuals equal to $k$,
and the number of $W$-colored individuals equal to $l$.
\smallskip

Let's denote by $Z(t)$ the Markov process  
with the state space $\mathcal L$ and the infinitesimal generator
$Q=(q_{\zeta \eta })$ given by the following formulas:
\begin{equation}\label{Q}
q_{\zeta \eta }=
\begin{cases}
\overline{x}\binom{k}{2},\quad&\text{if $\eta=(k-2,l+1)$, given $\zeta =(k,l)$}; \\ 
x\binom{k}{2}+xkl,\quad&\text{if $\eta =(k-1,l)$, given $\zeta =(k,l)$};\\ 
\overline{x}\binom{l}{2}+\overline{x}kl,\quad&\text{if $\eta =(k,l-1)$, given $\zeta =(k,l)$}; \\ 
x\binom{l}{2},\quad&\text{if $\eta =(k+1,l-2)$, given $\zeta=(k,l)$}; \\ 
-\binom{k+l}{2},\quad&\text{if $\eta =(k,l)$, given $\zeta =(k,l)$}; \\ 
0,&\text{otherwise.}
\end{cases}
\end{equation}
Here and in what follows, we take $\binom{n}{m}=0$ if $n<m$.

Notice that $Q$ is a square matrix with the size $n(n+3)/2$. We write it in 
the form of a block matrix with diagonal $(-r_{n}I_{n+1}, -r_{n-1}I_{n}, ..., -r_{2}I_{3}), 0_{2})$ and updiagonal $(B_{n+1,n}, B_{n,n-1}, ..., B_{3,2})$, where $r_{k}=\binom{k}{2},$ $I_{k}$ is the $k\times k$ identity matrix, 
$0_{2}$ is the $2\times2$ zero square matrix, and $B_{k,k-1}$ is a $k\times (k-1)$ matrix representing the transition (coalescent event) from the diagonal
$\Delta_{k-1}=\{(l_1,l_2)\in\mathcal L\,;\,l_1+l_2=k-1\}$
to the diagonal $\Delta_{k-2}=\{(l_1,l_2)\in\mathcal L\,;\,l_1+l_2=k-2\}$
of $\mathcal L$. All
other entries of $Q$ are zero.
\smallskip

Associated with the colored coalescent process $Z(t)$, there is
a discrete time Markov chain called the {\it jump chain}
of $Z(t)$. Its has the same state space $\mathcal L$ as that of
$Z(t)$ and its probability transition matrix is denoted by $J$.

The jump chain of $Z(t)$ has two absorbing states $(0,1)$ and $(1,0)$. 
When the coalescent process $Z(t)$ undergoes 
$k$ coalescent events, the associated
Markov chain moves $k$ steps. Thus, if we starts at
a state $(n_1,n_2)\in\mathcal L$ with $n_1+n_2=n$,
the probability of arriving at a state $(l_{1},l_{2})\in\mathcal L$ with 
$l_{1}+l_{2}=n-k$ can be read out from 
a row of $J^{k}$.
\smallskip

We denote
$$C_{n,n-k+1}:=r_{n}^{-1}r_{n-1}^{-1}\cdots r_{n-k+1}^{-1}B_{n+1,n}\cdots
B_{n,n-1}B_{n-k+2,n-k+1},$$
for $k\leq n-1.$ This is a $(n+1)\times (n-k+1)$ matrix. To fix this matrix,
we assume that the states are listed 
as $(0,n), (1,n-1), (2,n-2),\dots,(n,0),(0,n-1),(1,n-2), 
\dots,(0,1),(1,0)$. 

\begin{lemma} If $Z(t)$
starts at $(n_{1},n_{2})$, where $n_{1}+n_{2}=n$,
the probability that it reaches $(l_{1},l_{2})$, 
where $l_{1}+l_{2}=n-k$, after $k$ coalescent events is the 
$(n_{1}+1,l_{1}+1)$-entry of $C_{n,n-k+1}$. 
\end{lemma}

Notice that in particular, 
if $Z(t)$ starts at $(n_{1},n_{2})$, where \ $n_{1}+n_{2}=n$, the
probability that it reaches $(0,1)$ is the $
(n_{1}+1,1)$-entry of $C_{n,2}$, and the probability that it reaches $(1,0)$ 
is $(n_{1}+1,2)$-entry of $C_{n,2}.$
This means that after $n-1$ coalescent events, the process $Z(t)$ 
will be absorbed by either $(0,1)$ or $(1,0)$ and 
the absorption probabilities are given by the matrix 
$C_{n,2}$.

\begin{proof} By a direct computation, we see that the block of $J^k$
corresponding to the transition from $\Delta_n$ to $\Delta_{n-k}$ is 
exactly $C_{n,n-k+1}$. 
Then starting at a state $(n_{1},n_{2})$
with $n_1+n_2=n$, the $(n_{1}+1)$-th row of $C_{n,n-k+1}$ give us the
probability distribution of the process $Z(t)$ reaching
the state $(l_{1},l_{2})$ with $l_{1}+l_{2}=n-k$ after $k$ coalescent events.  
\end{proof}

We calculate next the expectation of the time $T_\pi$ for 
the coalescent process $Z(t)$ to reach an absorbing state, given that it
starts at an initial distribution 
$\pi$ on $\Delta_n$. For a non-absorption state $(k,l)\in\mathcal L$, let
$a_{\pi,(k,l)}$ be the the {\it sojourn coefficient} of $(k,l)$, which is  
the expected number of times the associated jump chain of $Z(t)$ 
to visit the state $(k,l)$, given that it starts at a
distribution $\pi$ on $\Delta_n$. Furthermore, let
$\omega_{(k,l)}$ be a random variable distributed exponentially
with the mean $r_{k+l}^{-1}$. The random variable $\omega_{(k,l)}$ 
is the time before the process $Z(t)$
to reach a state in $\Delta_{k+l-1}$, given that it starts at the state 
$(k,l)$. These random variables $\omega_{(k,l)}$, $(k,l)\in\mathcal L$
and $k+l\geq2$, are independent. Finally,
we define a random variable
\begin{equation}\label{coaltime}
\widetilde{T}_{\pi }=\sum_{m=2}^{n}\sum_{k+l=m}a_{\pi ,\left( k,l\right) }\omega _{k,l}
\end{equation}
\smallskip

\begin{lemma}\label{exp} The expectation of the coalescent time
is given by $$E(T_{\pi})=\sum_{m=2}^{n}\frac{1}{r_m}=2-\frac2{n}.$$
\end{lemma}

\begin{proof} The Feller relation (see \cite{sy}) says $E(T_\pi)=E(\widetilde{T}_\pi)$.
By a calculation involving the fundamental matrix for the associated 
Markov chain, it turns
out that $\sum_{k+l=m}a_{\pi ,\left( k,l\right)}=1$ for each $m\geq 2$.
(This matrix computation is straightforward but quite tedious. On the other hand, we shall see very easily while it is true in the proof of Theorem \ref{main1}.) Since 
$E(\omega _{k,l})=r_{k+l}^{-1}$, we get the desired expectation of 
$T_{\pi}.$
\end{proof}

\bigskip

\section{The lumpability of Markov processes}

The full details of the coalescent process $Z(t)$ is hard to compute
in general. For example, there are two absorption states, $(0,1)$ and $(1,0)$.
The coalescent time $T_\pi$ is the time elapsed before the process $Z(t)$
reaches any one of these two states. Can we compute the coalescent time to a particular one of these two states? In order to do this, the method of
lumping seems to be quite suitable. In what follows, we will show that $Z(t)$
can be lumped into another Markov process. For this lumped Markov process,
the computation of most of its interesting parameters can be made explicit. 
For that purpose,
we need to discuss the notion of lumpability of a Markov process first.
\smallskip

Let $X_{n}$ be a Markov chain with finite state space $S=\left\{
e_{1},e_{2,}\cdots \cdots ,e_{r}\right\}$, and $\overline{S}=$ $\left\{
E_{1},E_{2},\cdots \cdots ,E_{v}\right\} $ be a partition of $S$.
Let $p_{ij}$ be the transition probability from $e_i$ to $e_j$.
The probability that the chain moves into the set $E_\eta$ in one step,
given that it starts at $e_i$, is equal to $\sum_{e_{k}\in E_\eta}p_{ik}$. 
We say that the Markov chain $X_n$ is {\it lumpable} if 
for $e_i,e_j\in E_\xi$,
$$\sum_{e_{k}\in E_\eta}p_{ik}=\sum_{e_{k}\in E_\eta}p_{jk}.$$
When the chain $X_n$ is lumpable, we define a Markov chain $\overline{X}_n$
on $\overline{S}$ with the transition probability from $E_\xi$ to $E_\eta$
equal to $p_{\xi\eta}=\sum_{e_{k}\in E_\eta}p_{ik}$, for $e_i\in E_\xi$. This chain
$\overline{X}_n$ is  
called a {\it lumping} of $X_n$. See \cite{KS}.
\smallskip

Now, we consider a Markov process $X(t)$ on $S$. 
Let $P(t)=(p_{ij}(t))$ be the probability transition matrix 
of $X(t)$. We define that $X(t)$ is {\it lumpable} if
for $e_i,e_j\in E_\xi$,
$$\sum_{e_{k}\in E_\eta}p_{ik}(t)=\sum_{e_{k}\in E_\eta}p_{jk}(t),\quad
\text{for all $t\geq0$}.$$
\smallskip

Suppose that $X(t)$ is lumpable. We define
$$p_{\xi\eta}(t)=\sum_{e_k\in E_{\eta}}p_{ik}(t),$$
for $e_i\in E_\xi$. Let $\overline{P}(t)=(p_{\xi\eta}(t))$.
This is an $v\times v$ matrix.
\smallskip

\begin{lemma} $\overline{P}(t)$ defines a Markov process on $\overline{S}$.
\end{lemma}

\smallskip

The proof is to check the property that $\overline{P}(t)\overline{P}(t+s)$.  It can be done directly by computation.  We call the Markov process $\overline{X}(t)$ on $\overline{S}$ with the probability
transition matrix $\overline{P}(t)$ a {\it lumping} of $X(t)$.
\smallskip

Let $Q=(q_{ik})$ be the infinitesimal generator of $X(t)$, i.e. $P(t)=e^{tQ}$. 
We say that $Q$ is {\it lumpable} if 
$$\sum_{e_{k}\in E_{\eta}}q_{ik}=\sum_{e_{k}\in E_{\eta}}q_{jk}$$
for $e_i,e_j\in E_\xi$.
\smallskip

We need to introduce some notations in order to prove the
following theorem which relates the lumpability of $Q$ with that
of $P(t)=e^{tQ}$. Recall that the finite set of states
$S=\{e_1,e_2,\dots,e_r\}$ is partitioned into $\overline{S}=\{
E_1,E_2,\dots,E_v\}$. Let $U$ be the $v\times r$ matrix
whose $\xi$-th row is the probability vector having equal components for
states in $E_\xi$, $\xi=1,2\dots,v$, and $0$ elsewhere. Let $V$ be the $r\times
v $ matrix with the $\eta$-th column a vector with 1's in the components
corresponding to states in $E_\eta$, $\eta=1,2,\dots,v$ and $0$ elsewhere. 
One can check by a direct computation that $UV=I_v$ and that the Markov process $X(t)$ is 
lumpable if and only if $VUP(t)V=P(t)V$. (This statement generalizes a theorem about lumpability of Markov chains in \cite{KS}.)
Furthermore, when $P(t)$ is lumpable,
we have the $\overline{P}(t)=UP(t)V$ for the lumping process $\overline{X}(t)$. 
Similarly, the infinitesimal generator $Q$ of $X(t)$ is lumpable if and only if $VUQV=QV$. 
And when $Q$ is lumpable, its lumping is $\overline{Q}=UQV$.

\begin{theorem}\label{lumpable}
A necessary and sufficient condition for the Markov process $X(t)$ to be lumpable is that
its infinitesimal generator $Q$ is lumpable. When $Q$ is lumpable, we have
$\overline{P}(t)=e^{t\overline{Q}}$.

\end{theorem}

\smallskip 
The intuitive idea of this theorem is clear, we therefore will not give the detailed proof here.

\section{Parity lumping of the colored coalescent process}

The state space $\mathcal L$ of $Z(t)$ is partitioned into diagonals
$\Delta_m$, $m=1,2,\dots,n$. We will divide each $\Delta_m$ into two
disjoint subsets, $O_m$ and $E_m$. A state $(k,l)\in O_m$ when
$k+l=m$ and $k$ is odd and $(k,l)\in E_m$ when $k+l=m$ and $k$ is even.
Let
$$\overline{\mathcal L}=\{O_m,E_m\,;\,m=1,2,\dots,n\}.$$
We will define a new Markov process on $\overline{\mathcal L}$ obtained by  
a lumping of the process $Z(t)$. For that purpose, we need to
check the lumpability of $Z(t)$. Using Theorem \ref{lumpable},
we only need to check the lumpability of the infinitesimal
generator $Q=(q_{\zeta\eta})$ given by Equation (\ref{Q}). 
\smallskip

Let $\zeta=(k,l)\in O_m$, then $k+l=m$ and $k$ is odd. We have
$$\sum_{\eta\in O_{m-1}}q_{\zeta\eta}=\overline{x}\binom{k}{2}+\overline{x}
\binom{l}{2}+\overline{x}kl=\overline{x}\binom{m}{2}$$
and
$$\sum_{\eta\in E_{m-1}}q_{\zeta\eta}=x\binom{k}{2}+xkl+x\binom{l}{2}=x\binom{m}{2}.$$
Both of them are independent of $k$. Other lumpability conditions
for $Q$ can be checked similarly. So, by Theorem \ref{lumpable},
$Z(t)$ has a lumping $\overline{Z}(t)$, which is a Markov process
on $\overline{\mathcal L}$ with the infinitesimal generator $\overline{Q}$.
We will order the elements in $\overline{\mathcal L}$ by $E_n,O_n,E_{n-1},O_{n-1},\dots,
E_2,O_2,E_1,O_1$. Under this ordering, $\overline{Q}$ is a block matrix with diagonal $(-r_{n}I_2, ..., -r_{2}I_2, 0_2)$ and updiagonal $(r_{n}C, ..., r_{2}C)$, where $C$ is a $2\times 2$ matrix: 

$$C=\left[
\begin{array}{cc}
\overline{x}& x\\
x&\overline{x}
\end{array}\right].$$

We will follow the same idea as in Section 2 to compute 
basic parameters of the lumped coalescent process $\overline{Z}(t)$.
We shall be able to get more complete information about the process
$\overline{Z}(t)$ than about $Z(t)$.
\smallskip

\noindent(1) {\it The fundamental matrix of the jump chain 
of $\overline{Z}(t)$.} The probability transition matrix of the jump chain of $\overline{Z}(t)$
is a block matrix with each updiagonal entry $C$ and each diagonal $0$ except the last $I_2$.
Notice that the states $E_1$ and $O_1$ are absorbing states. So the fundamental matrix 
$\overline{N}$ of the jump chain can be easily calculated.  This is a block updiagonal matrix, each of whose diagonal entry is $I_2$, each of whose updiagonal entry is $C$, each of whose second updiagonal entry is $C^2$, up to the entry  $C^{n-2}$.  By a straightforward computation, we may get 

$
C^k=\left[
\begin{array}{cc}
\frac{1}{2}+\frac{1}{2}(1-2x)^{k} & \frac{1}{2}-\frac{1}{2}(1-2x)^{k} \\ 
\frac{1}{2}-\frac{1}{2}(1-2x)^{k} & \frac{1}{2}+\frac{1}{2}(1-2x)^{k}
\end{array}
\right] 
$
\smallskip

\noindent(2) {\it The sojourn coefficients of the jump chain.} The sojourn
coefficients $a_k$ (respectively, $b_k$) is
the expected number of times the jump chain $\overline{J}$ of $\overline{Z}(t)$ 
to visit the state $E_k$ (respectively, $O_k$), given that it starts at an
initial distribution $\pi=(\pi_E,\pi_O)$ on the states $\{E_n,O_n\}$. 
\smallskip

\begin{lemma} We have
\begin{equation}
(a_{k},b_{k})=\left( \frac{1}{2}+\frac{1}{2}\left( \pi _{E}-\pi _{O}\right) (1-2x)^{n-k},\,
\frac{1}{2}+\frac{1}{2}\left( \pi _{O}-\pi _{E}\right) (1-2x)^{n-k}\right) 
\end{equation}
\end{lemma}
\smallskip

\begin{proof} The sojourn coefficients $(a_k,b_k)$ can be calculated as follows:
$$\begin{aligned}
(a_{k},b_{k})&=\left( \pi _{E},\pi _{O}\right) C^{n-k}\\
&=\left( \pi _{E},\pi _{O}\right) \left( 
\begin{array}{cc}
\frac{1}{2}+\frac{1}{2}(1-2x)^{n-k} & \frac{1}{2}-\frac{1}{2}(1-2x)^{n-k} \\ 
\frac{1}{2}-\frac{1}{2}(1-2x)^{n-k} & \frac{1}{2}+\frac{1}{2}(1-2x)^{n-k}
\end{array}
\right) \\
&=\left( \frac{1}{2}+\frac{1}{2}\left( \pi _{E}-\pi _{O}\right) (1-2x)^{n-k},
\frac{1}{2}+\frac{1}{2}\left( \pi _{O}-\pi _{E}\right) (1-2x)^{n-k}\right).
\end{aligned}$$
\end{proof}
\smallskip

\noindent(3) {\it Coalescent probability for the lumped coalescent process $\overline{Z}(t)$.}
Since $E_1$ and $O_1$ are absorbing states, the coalescent probability $P_{\pi,E}$ (respectively, $P_{\pi,O}$) 
is the probability of the jump chain to arrive at $E_1$ (respectively, $O_1$), given that it starts at a
distribution $\pi=(\pi_E,\pi_O)$ on the initial states $\{E_n,O_n\}$.  
\smallskip

\begin{lemma} We have 
$$
(P_{E},P_{O})=\left( \frac{1}{2}+\frac{1}{2}\left( \pi _{E}-\pi _{O}\right) (1-2x)^{n-1},\,
\frac{1}{2}+\frac{1}{2}\left( \pi _{O}-\pi _{E}\right) (1-2x)^{n-1}\right). 
$$
In particular, let $P_{E,E}$ be the probability of reaching $E_1$, given that the process starts at $E_n$,
and other quantities $P_{E,O},P_{O,E},P_{O,O}$ be defined similarly, then we have
\begin{equation}\label{abP}
\begin{aligned}
&P_{E,E}=\frac{1}{2}+\frac{1}{2}(1-2x)^{n-1}, & P_{E,O}=\frac{1}{2}-\frac{1}{2}(1-2x)^{n-1},\\
&P_{O,E}=\frac{1}{2}-\frac{1}{2}(1-2x)^{n-1}, &P_{O,O}=%
\frac{1}{2}+\frac{1}{2}(1-2x)^{n-1}
\end{aligned}
\end{equation}
\end{lemma}
\smallskip

\begin{proof} Since a minimal Markov process on a finite set of states and its jump chain have the same character of states, $E_1$ and $O_1$ are absorbing states for both $\overline{Z}(t)$ and $\overline{J}$. So we have  $(P_{\pi,E},P_{\pi,O})=(a_1,b_1)$.
\end{proof}
\smallskip

\noindent(4) {\it The expected coalescent time for the lumped coalescent process $\overline{Z}(t)$.} Since $E_{1}$ and $O_{1}$ both are absorbing states
of the process $\overline{Z}(t)$, the coalescent time $T_{\pi,E}$ to the state $E_{1}$ is a random variable, which is the time elapsed before the process $\overline{Z}(t)$ reaches $E_{1}$, given that it starts at the initial distribution $\pi=(\pi_E,\pi_O)$ and conditional on $\overline{Z}(t)$ not reaching $O_{1}$. Similarly, we have the coalescent time $T_{\pi,O}$ to the state $O_1$ as a random variable. To calculate the expectation of the coalescent time, we will define a random variable $\widetilde{T}_{\pi ,E}$ first and calculate its expectation. Then the Feller relation will tell us that the expectation of $\widetilde{T}_{\pi,E}$ is exactly the same as the expectation of the coalescent time $T_{\pi,E}$ to the states $E_1$. I. e., we have $E(\widetilde{T}_{\pi,E})=E(T_{\pi,E})$.

Let 
$\tau _{k}$ be a random variable distributed exponentially with the mean $r_{k}^{-1}$, which is the time elapsed before the lumped coalescent process 
$\overline{Z}(t)$ moves either to $E_{k-1}$ or $O_{k-1}$, given that it starts at $E_{k}$ (or $O_{k}$). Recall that the sojourn coefficients $a_{k}$
(respectively, $b_{k}$) is the expected number of times the jump chain $\overline{J}$ of $\overline{Z}(t)$ to visit the state $E_{k}$ (respectively, 
$O_{k}$), given that it starts at a initial distribution $\pi =(\pi _{E},\pi
_{O})$ on the initial states \{$E_{n},O_{n}$\}. Then the random variable $\widetilde{T}_{\pi,E}$ is defined by
\begin{equation}
\begin{aligned}
\widetilde{T}_{\pi,E}&=\sum_{k=3}^{n}a_{k}\tau _{k}+\sum_{k=3}^{n}b_{k}\tau _{k}+
\frac{a_{2}}{\overline{x}}\tau _{2}+\frac{b_{2}}{x}\tau _{2} \\
&=\sum_{k=3}^{n}\tau _{k}+\frac{1}{2\overline{x}x}\left(1+(\pi _{O}-\pi
_{E})(1-2x)^{n-1}\right)\tau _{2}.
\end{aligned}
\end{equation}
\smallskip

Similarly, we may define a random variable $\widetilde{T}_{\pi ,O}$, whose expectation turns out to be the same as the expectation of the coalescent time $T_{\pi,O}$ to the state $O_1$: 
\begin{equation}
\widetilde{T}_{\pi,O}=\sum_{k=3}^{n}\tau _{k}+\frac{1}{2\overline{x}x}\left(1+(\pi _{E}-\pi
_{O})(1-2x)^{n-1}\right)\tau _{2}.
\end{equation}
\smallskip

\begin{lemma} The expectation of the coalescent time to $E_{1}$ and $O_{1}$
are 
\begin{equation*}
E(T_{\pi ,E})=1-\frac{2}{n}+\frac{1}{2\overline{x}x}\left(1+(\pi _{O}-\pi
_{E})(1-2x)^{n-1}\right)
\end{equation*}
and
\begin{equation*}
E(T_{\pi ,O})=1-\frac{2}{n}+\frac{1}{2\overline{x}x}\left(1+(\pi _{E}-\pi
_{O})(1-2x)^{n-1}\right)
\end{equation*}
respectively. 
\end{lemma}

\begin{proof}
The expectations 
$$
\begin{aligned}
&E(\widetilde{T}_{\pi ,E})=1-\frac{2}{n}+\frac{1}{2\overline{x}x}\left(1+(\pi _{O}-\pi
_{E})(1-2x)^{n-1}\right)\\
&E(\widetilde{T}_{\pi ,O})=1-\frac{2}{n}+\frac{1}{2\overline{x}x}\left(1+(\pi _{E}-\pi
_{O})(1-2x)^{n-1}\right)
\end{aligned}
$$ are easy to compute since we know that $E(\tau _{k})=r_{k}^{-1}$ and 
$\tau _{k}$'s are independent random variables. To see why 
$$E(\widetilde{T}_{\pi ,E})=E(T_{\pi,E}),$$ 
we define a
conditional process whose state characters are all the same as that of the
process $\overline{Z}(t)$ except for the states $E_{2}$ and $O_{2}$. At these two states $E_2$ and $O_2$, the conditional process will move only to $E_1$ with the mean of holding time $\frac{1}{\overline{x}r_{2}}$ and $\frac{1}{xr_{2}}$,
respectively. Then, by the Feller relation, $E(\widetilde{T}_{\pi,E})$ is the expectation of the coalescent time for this conditional process to 
reach $E_{1}$. On the other hand, the coalescent time for the conditional process to reach $E_1$ is exactly $T_{\pi,E}$. So we have $E(\widetilde{T}_{\pi,E})=E(T_{\pi,E})$. Similarly, we have $E(\widetilde{T}_{\pi,O})=E(T_{\pi,O})$.
\end{proof}

\noindent(5) {\it The cumulative distribution function of the coalescent time of
$\overline{Z}(t)$.} The computation of the complementary cumulative distribution functions
$\Pr\,\{T_{\pi,E}\geq t\}$ and $\Pr\,\{T_{\pi,O}\geq t\}$ is rather involved. So we will put it in the appendix and give here only the results.
\smallskip

We need to introduce some functions of $x$ first. They are $K_{n,2}$, $K_{n,2'}$, $K_{n',2}$, and $K_{n',2'}$. Here the indices $1,2,\dots,n$ are for the states
$E_1,E_2,\dots,E_n$ and the indices $1',2',\dots,n'$ are for $O_1,O_2,\dots,O_n$. These functions are given as follows:
\begin{equation}
\begin{aligned}
&K_{n,2}=\left(\frac{1}{2}+\frac{1}{2}(1-2x)^{n-2}\right)\sum_{k=0}^{n-3}\frac{r_{n-k}}{
r_{n-k}-\overline{x}r_{2}}\prod
_{\begin{subarray}{l}i=0\\ i\neq k\end{subarray}}^{n-3}\frac{r_{n-i}}{r_{n-i}-r_{n-k}},  \\
&K_{n,2'}=\left(\frac{1}{2}-\frac{1}{2}(1-2x)^{n-2}\right)\sum_{k=0}^{n-3}\frac{r_{n-k}
}{r_{n-k}-xr_{2}}
\prod_{\begin{subarray}{l}i=0\\ i\neq k\end{subarray}}^{n-3}\frac{r_{n-i}}{r_{n-i}-r_{n-k}},\\
&K_{n',2}=\left(\frac{1}{2}-\frac{1}{2}(1-2x)^{n-2}\right)\sum_{k=0}^{n-3}\frac{r_{n-k}%
}{r_{n-k}-\overline{x}r_{2}}
\prod_{\begin{subarray}{l}i=0\\ i\neq k\end{subarray}}^{n-3}\frac{r_{n-i}}{r_{n-i}-r_{n-k}},\\
&K_{n',2'}=\left(\frac{1}{2}+\frac{1}{2}(1-2x)^{n-2}\right)\sum_{k=0}^{n-3}\frac{r_{n-k}}{r_{n-k}-xr_{2}}\prod_{\begin{subarray}{l}i=0\\ i\neq k\end{subarray}}^{n-3}\frac{r_{n-i}}{r_{n-i}-r_{n-k}}.
\end{aligned}
\end{equation}
\smallskip

We denote $T_{\pi,E}$ with $\pi_E=1,\pi_O=0$ by $T_{E,E}$. Other coalescent time $T_{O,E}$, $T_{E,O}$, and $T_{O,O}$ have the same meaning. 
\smallskip

\begin{lemma}\label{hard} We have
\begin{equation*}
\begin{aligned}
\Pr&\,\{T_{E,E}\geq t\}= \\
&\sum_{k=0}^{n-3}\,\prod_{\begin{subarray}{l}i=0\\ i\neq k\end{subarray}}^{n-3}\frac{r_{n-i}}{r_{n-i}-r_{n-k}}\left(
\frac{\overline{x}(\frac{1}{2}+\frac{1}{2}(1-2x)^{n-2})}{\overline{x}
r_{2}-r_{n-k}}+\frac{x(\frac{1}{2}-\frac{1}{2}(1-2x)^{n-2})}{xr_{2}-r_{n-k}}%
\right)e^{-r_{n-k}t} \\
&\qquad+K_{n,2}e^{-\overline{x}r_{2}t}+K_{n,2'}e^{-xr_{2}t}
\end{aligned}
\end{equation*}
and
\begin{equation*}
\begin{aligned}
\Pr&\,\{T_{O,E}\geq t\}= \\
&\sum_{k=0}^{n-3}\,\prod_{\begin{subarray}{l}i=0\\ i\neq k\end{subarray}}^{n-3}\frac{r_{n-i}}{r_{n-i}-r_{n-k}}\left(
\frac{\overline{x}(\frac{1}{2}+\frac{1}{2}(1-2x)^{n-2})}{\overline{x}
r_{2}-r_{n-k}}+\frac{x(\frac{1}{2}-\frac{1}{2}(1-2x)^{n-2})}{xr_{2}-r_{n-k}}
\right)e^{-r_{n-k}t} \\
&\qquad+K_{n',2}e^{-\overline{x}r_{2}t}+K_{n',2'}e^{-xr_{2}t}.
\end{aligned}
\end{equation*}

Also, $\Pr\,\{T_{E,O}\geq t\}=\Pr\,\{T_{O,E}\geq t\}$ and $\Pr\,\{T_{O,O}\geq t\}=\Pr\,\{T_{E,E}\geq t\}$.
\end{lemma}
\smallskip

\section{Back to the colored coalescent process $Z(t)$}

For the parity lumping of the colored coalescent process, we have a
commutative diagram 
\begin{equation*}
\CD
Z(t)  @>\text{jump chain}>>  J \\ 
@V\text{lumping}VV  @VV\text{lumping}V\\ 
\overline{Z}(t) @>\text{jump chain} >> \overline{J}.
\endCD
\end{equation*}
This can be shown by a simple
matrix computation. Notice that, for an arbitrary lumpable Markov process, we do not know
if this diagram is commutative or under what condition this diagram is
commutative. 
\smallskip

The commutativity of the above diagram for the parity lumping of our colored coalescent process $Z(t)$ provides a way for us to recover information about $Z(t)$ from our knowledge of the lumped coalescent process $\overline{Z}(t)$.
Due to the fact that both processes $Z(t)$ and $\overline{Z}(t)$ have the same
absorbing states, we can achieve a complete recovery of information for certain parameters of the coalescent process $Z(t)$.
\smallskip

Let $N$ be the fundamental matrix of the jump chain $J,$ of $Z(t)$, and 
$\overline{N}$ be the fundamental matrix of the jump chain $\overline{J}$ of $\overline{Z}(t)$. Notice that the absorption states of $\overline{J}$ is $E_1=\{(0,1)\}$ and $O_1=\{(1,0)\}$.
Let $U_{0}$ and $V_{0}$ be
matrices obtained from $U$ and $V$ by dropping the last two rows and columns
corresponding to the states $(0,1)$ and $(1,0)$. Then it is easy to verify that
$N$ is lumpable by $U_0$ and $V_0$: $V_0U_0NV_0=NV_0$. Furthermore, we have
$\overline{N}=U_{0}NV_{0}$. 
\smallskip

The following is the main theorem of this section.

\begin{theorem}\label{main1}
Let $P_{(n_1,n_2)(0,1)}$ and $P_{(n_1,n_2)(1,0)}$ be the probabilities that the coalescent process $Z(t)$ reaches $(0,1)$ and $(1,0)$, respectively,
given that it starts at $(n_1,n_2)$, $n_1+n_2=n$. Then we have
\begin{equation}
P_{(n_{1},n_{2})(0,1)}=\begin{cases}
\frac{1}{2}+\frac{1}{2}(1-2x)^{n-1},&\text{if $n_{1}$ is even,} \\ 
\frac{1}{2}-\frac{1}{2}(1-2x)^{n-1}, &\text{if $n_{1}$ is odd;} 
\end{cases}
\end{equation}
\begin{equation}
P_{(n_{1},n_{2})(1,0)}=\begin{cases}
\frac{1}{2}-\frac{1}{2}(1-2x)^{n-1},&\text{if $n_{1}$ is even,} \\ 
\frac{1}{2}+\frac{1}{2}(1-2x)^{n-1},&\text{if $n_{1}$ is odd. }
\end{cases}
\end{equation}

Furthermore, let $T_{(n_1,n_2)(0,1)}$ (or $T_{(n_1,n_2)(1,0)}$) be the 
time for the coalescent process $Z(t)$ to reach
$(0,1)$ (or $(1,0)$), given that it starts at $(n_1,n_2)$, $n_1+n_2=n$. 
Then the expectations of $T_{(n_1,n_2)(0,1)}$ and $T_{(n_1,n_2)(1,0)}$
are given as follows:
\begin{equation}
E(T_{(n_1,n_2)(0,1)})=\begin{cases}
1-\frac{2}{n}+\frac{1}{2\overline{x}x}\left(1-(1-2x)^{n-1}\right),&\text{if $n_1$ is even,}\\
1-\frac{2}{n}+\frac{1}{2\overline{x}x}\left(1+(1-2x)^{n-1}\right),&\text{if $n_1$ is odd;}
\end{cases}
\end{equation}

\begin{equation}
E(T_{(n_1,n_2)(1,0)})=\begin{cases}
1-\frac{2}{n}+\frac{1}{2\overline{x}x}\left(1+(1-2x)^{n-1}\right),&\text{if $n_1$ is even,}\\
1-\frac{2}{n}+\frac{1}{2\overline{x}x}\left(1-(1-2x)^{n-1}\right),&\text{if $n_1$ is odd.}
\end{cases}
\end{equation}
\end{theorem}

\begin{proof} 
The absorption probabilities of the jump chains $J$ and $\overline{J}$ are  calculated from the matrices
$$N\left[\begin{array}{c}
0\\ 
\vdots\\
0\\
C
\end{array}\right]\quad\text{and}\quad
\overline{N}
\left[\begin{array}{c}
0\\ 
\vdots\\
0\\
C
\end{array}\right]
=U_0NV_0\left[\begin{array}{c}
0\\ 
\vdots\\
0\\
C
\end{array}\right]=U_0N\left[\begin{array}{c}
0\\ 
\vdots\\
0\\
C
\end{array}\right],
$$
respectively. So $P_{(n_1,n_2)(0,1)}$, for $n_1+n_2=n$ and $n_1$ even, is equal to $P_{E,E}$.
By Equation (\ref{abP}), we get the desired value for $P_{(n_1,n_2)(0,1)}$ in this
case. All other cases can be obtained similarly.
\smallskip

To compute the coalescent time, we first denote by $a_{(n_1,n_2)(k_1,k_2)}$
the sojourn coefficient of the jump chain $J$, which is the expect number of times the jump chain $J$ visits the state $(k_1,k_2)$, given that it
starts at the state $(n_1,n_2)$, $n_1+n_2=n$. Since the jump chain $J$ is lumpable and its lumping is the jump chain $\overline{J}$, we have
\begin{equation}\label{coef}
\begin{aligned}
&\sum_{k_1+k_2=k,\,k_1\text{\,even}}a_{(n_1,n_2)(k_1,k_2)}=a_k\\
&\sum_{k_1+k_2=k,\,k_1\text{\,odd}}a_{(n_1,n_2)(k_1,k_2)}=b_k,
\end{aligned}
\end{equation}
where $a_k$ and $b_k$ are the sojourn coefficients of $\overline{J}$ 
corresponding to the expected number of times $\overline{J}$ visits 
$E_k$ and $O_k$, respectively, given that it starts at the distribution
\begin{equation}\label{pi}
\pi=(\pi_E,\pi_O)=\begin{cases}
(1,0),&\text{if $n_1$ is even,}\\
(0,1),&\text{if $n_1$ is odd.}
\end{cases}
\end{equation}
Thus, we have
$$\begin{aligned}
\widetilde{T}_{(n_1,n_2)(0,1)}&=\sum_{k=3}^n\left(\,\sum_{k_1+k_2=k}a_{(n_1,n_2)(k_1,k_2)}\right)\tau_k\\ &\qquad+\frac{a_{(n_1,n_2)(0,2)}}{\overline{x}}\tau_2+\frac{a_{(n_1,n_2)(2,0)}}{\overline{x}}\tau_2
+\frac{a_{(n_1,n_2)(1,1)}}{x}\tau_2\\
&=\sum_{k=3}^na_k\tau_k+\sum_{k=3}^nb_k\tau_k+\frac{a_2}{\overline{x}}\tau_2+\frac{b_2}{x}\tau_2\\
&=\widetilde{T}_{\pi,E},
\end{aligned}
$$ 
where $\pi$ is the distribution given by Equation (\ref{pi}). So Lemma 6 gives us the desired expectation of $T_{(n_1,n_2)(0,1)}$. The calculation of
$E(T_{(n_1,n_2)(1,0)})$ is exactly the same.
\end{proof}
\smallskip

\begin{remark} Equation (\ref{coef}) shows $\sum_{k+l=m}a_{\pi ,\left( k,l\right)}=1$ as claimed in the proof of Lemma \ref{exp}.
\end{remark}
\smallskip

\begin{remark} For $(n_1,n_2)\in\Delta_n$ with $n_1$ even, it is obvious that $T_{(n_1,n_2)(0,1)}=T_{E,E}$. Thus Lemma 2.7 also gives us the complementary cumulative distribution function of $T_{(n_1,n_2)(0,1)}$. All other cases are the same so that the cumulative distribution functions of $T{(n_1,n_2)(0,1)}$ and $T_{(n_1,n_2)(1,0)}$ are all known.  
\end{remark}
\smallskip

There is, of course, information of the coalescent process $Z(t)$ which can not be recovered from what we know about $\overline{Z}(t)$. 
For example, recall that the parity of a state $(k,l)\in\mathcal L$ is the
parity (even or odd) of the integer $k$. Let us define $\rho_{k}$ to be the parity of a state after the $k$-th coalescent
event. In other words, $\rho_{k}$ is a random variable which takes the value 0 if the coalescent process $Z(t)$ is in $E_k$ and takes the value 1 if the coalescent process is in $O_k$, after the $k$-th coalescent event. We will
use $\rho_0$ to denote the random variable which takes the value 0 if the parity
of the initial state is even and the value 1 otherwise. Then we have
\begin{equation}\label{parity}
\begin{aligned}
&\text{Pr}(\rho_k=0\,|\,\rho_0=0)=\frac{1}{2}+\frac{1}{2}(1-2x)^{k}\\
&\text{Pr}(\rho_k=1\,|\,\rho_0=0)=\frac{1}{2}-\frac{1}{2}(1-2x)^{k}\\
&\text{Pr}(\rho_k=0\,|\,\rho_0=1)=\frac{1}{2}-\frac{1}{2}(1-2x)^{k}\\ 
&\text{Pr}(\rho_k=1\,|\,\rho_0=1)=\frac{1}{2}+\frac{1}{2}(1-2x)^{k}
\end{aligned}
\end{equation}
\smallskip

From Equation (\ref{parity}), we can get the probability of sequences of $n-1$ coalescent events with the parities of the passing states all specified. For example, the probability that all of $n-1$ coalescent events happen in states
with even parity, given that the parity of the initial state is also even, is
$$\text{Pr}(\rho_k=0\,;\,k=0,1,2,\cdots ,n-1)=\prod_{k=0}^{n-1}\left( \frac{1}{2}+\frac{1}{2}\left( 1-2x\right) ^{k}\right).$$
The probability of other sequences of coalescent events with specified parities 
can be calculated similarly.
\smallskip

On the other hand, starting at $(n_{1},n_{2})$, where \ $n_{1}+n_{2}=n$, after $k$ coalescent events, we do not know exactly the distribution of the
states $(l_{1},l_{2})$, $l_{1}+l_{2}=n-k$. The knowledge of the distribution of the parities of these states is not enough for recovering the distribution of these states.\\ 

\bigskip
\noindent {\it Acknowledgments.}  We would like to thank Prof. Michael Clegg for introducing us to the coalescent theory and its significance in population genetics.  The first author also would like to thank Prof. Bai-Lian Li for his encouragement.  We also thank the referees for their comments.

\medskip

\end{document}